\documentclass{article}
\usepackage{maa-monthly}
\usepackage{ulem}
\usepackage{hyperref}
\usepackage{nag}
\usepackage{color}

\theoremstyle{theorem}

\theoremstyle{definition}

\begin{document}

\title{Skip Letters for Short Supersequence of All Permutations}
\markright{Skip Letters}
\author{Oliver Tan}

\maketitle

\begin{abstract}
A supersequence over a finite set is a sequence that contains as subsequence all permutations of the set. This paper defines an infinite array of methods to create supersequences of decreasing lengths. This yields the shortest known supersequences over larger sets. It also provides the best results asymptotically. It is based on a general proof using a new property called strong completeness. The same technique also can be used to prove existing supersequences which combines the old and new ones into an unified conceptual framework.
\end{abstract}

\noindent
A subsequence $\sigma_1$ of a sequence $\sigma$ is a sequence that can be obtained from $\sigma$ by deleting none or some elements from $\sigma$, but leaving the order of the remaining elements intact. The relation will be denoted by $\sigma_1 < \sigma$. A sequence is called a supersequence over a finite set $A=\{a_1, a_2, ... a_n\}$ if it contains as subsequence all permutations of $a_1, a_2, ... a_n$. For examples, $\langle 1,2,3,1,2,1,3 \rangle$ and $\langle1,2,3,4,1,2,3,1,4,2,1,3 \rangle$ are supersequence over $\{1,2,3\}$ and $\{1,2,3,4\}$ respectively. They also happen to be the shortest supersequence over their sets.  A member of $A$ like $a_1$ is called a letter. 

Finding the shortest supersequence over a finite set is part of the larger universal permutation research area (Engen-Vatter [2]), related also to the shortest common supersequences and longest common subsequences issue with wide practical applications in biology, data compression, text editing and many others. It has the longest history among the various classes of problems in the area. Historically in the 1970s to early 1980, Newey [7], Adleman [1], Koutas-Hu [5], Galbiati-Preparata [3] and Mohanty [6] provide different algorithms to create supersequences of length $m^2 -2m + 4$ over a set of $m$ letters. For a long time, it was thought that this was the shortest that can be achieved (Koutas-Hu [5]). It stood to be the best result for over thirty years until Zalinescu [9] improves it by $1$ to the length of $m^2 -2m + 3$. Shortly, Radomirovic [8] further reduces the length to $\lceil m^2 -  \frac {7}{3}m + \frac{19}{3} \rceil$. This paper provides an alternative proof that enables generalization to remove unbounded many skip letters. The result is a shorter length that asymptotically approaching $\lceil m^2 - \frac {5}{2}m + C_{\epsilon} \rceil$. Formally, for any real number $\epsilon > 0$, there exists a constant $C_{\epsilon}$ such that for large enough $m$, there is a supersequence of length $\lceil m^2 - (\frac {5}{2}-\epsilon)m + C_{\epsilon} \rceil$ over a set of $m$ letters. With appropriate choice for the number of skip letters for each $m$, it is proved that there exists a supersequence of length $\lceil m^2 -  \frac {5}{2}m + \frac{3}{2}  (\frac{m}{2})^{\frac {2}{3}} +  (\frac{m}{2}) ^{\frac {1}{3}} + 7\rceil$.

It is still an open question what is the optimal length of a supersequence over a set of $m$ letters. On the lower bound side, Kleitman-Kwiatkowski [4] proves that a supersequence must have a minimum length of $m^2- C_{\epsilon} m^{7/4+\epsilon}$ for any $\epsilon > 0$, where $C_{\epsilon}$ is a constant depending on $\epsilon$.

\section{Notations and Basic Results}
\begin{flushleft} We begin with the notations and introduce the concept of strong completeness that can be applied to prove existing supersequences over finite sets. When we say $\sigma$ is a sequence over a set $A$, it means that all elements of $\sigma$ are letters from the set $A$. If $\sigma$ is a sequence, then $l_\sigma$ denotes the length of the sequence. For any $i$ where $1 \le i \le l_\sigma$, the notation $\sigma[i]$ denotes the $i$th element of $\sigma$, and $a = \sigma[i]$ is called an element of $\sigma$. An element $a$ of $\sigma$ is denoted by $a \in \sigma$. The position $i$, where the element $a$ occurs in the sequence $\sigma$, is denoted by $\sigma^{-1}[a]$. This definition is valid if $a$ occurs only once in the sequence $\sigma$. So when the definition is valid, we always have $\sigma[\sigma^{-1}[a]] = a$. The last element, second last element and so on of $\sigma$ are denoted by $\sigma[-1]$, $\sigma[-2]$ and so on, respectively. Given two integers $i$ and $j$, let $\sigma[i,j]$ represent the substring of $\sigma$ starting at position $i$ and ending at $j$, inclusively for both ends. For example, if $\sigma = \langle 1,2,3,4,5,6 \rangle$, then $\sigma[3,-2] = \langle 3,4,5\rangle$ represents the substring from the third element to the second last element of $\sigma$. For any $i$ where $1 \le i \le l_\sigma$, the shorter notation $\sigma | i$ is also used to represent $\sigma[1,i]$ to denote the substring of $\sigma$ starting from $1$ and ending at $i$ inclusive at both ends.  \end{flushleft}

Given a list of sequences $\sigma_1, \sigma_2, \sigma_3,...$, and letters $a_1, a_2, a_3,...$, then $\sigma_1 \sigma_2 \sigma_3...$ will denote the concatenation of those sequences $\sigma_1, \sigma_2, \sigma_3,...$, and $\sigma_1 a_1 \sigma_2 a_2 \sigma_3 a_3...$ will denote the concatenation of those sequences $\sigma_1, \sigma_2, \sigma_3,...$ interposed with those letters $a_1, a_2, a_3,...$. We may optionally write a dot $\cdot$ in between sequence and element, or between sequence and sequence, if it makes reading easier, like $\sigma_1 \cdot a_1  \cdot \sigma_2  \cdot a_2 \cdot \sigma_3 \cdot a_3...$.

For a set $A$, we use $|A|$ to denote the number of letters in $A$, $[A]_k$ to represent the set of all sequences over A of length $k$, where all elements of each sequence are distinct. We use $[A]$ to denote $[A]_n$ where $n = |A|$, i.e. $[A]$ is the set of all permutations over $A$. A sequence $\sigma$ is said to be  \textit{$k$-complete} for some $k \le |A|$ if for all $\sigma_1 \in [A]_k, \sigma_1 < \sigma$, i.e. all sequence of length $k$ containing only distinct elements is a subsequence of $\sigma$. A $|A|$-complete sequence is therefore a supersequence over $A$. 

\begin{flushleft} \textbf{Definition 1.} Suppose $\sigma_1, \sigma_2, ..., \sigma_n$ is a list of sequences over a set $A$, where $n \le |A|$. The list is said to be \textit{forward complete} if for any integer $k$ such that $1\le k \le n$, $\sigma_1 \sigma_2... \sigma_k$ is $k$-complete.  The list is said to be \textit{backward complete} if for any such $k$, $\sigma_{n-k+1} \sigma_{n-k+2}... \sigma_n$ is $k$-complete. The list is said to be \textit{strongly complete} if it is both forward complete and backward complete.
\end{flushleft}

An obvious corollary is that for any forward or backward complete list, $\sigma_1 \sigma_2, ... \sigma_n$ is $n$-complete. The simplest example of strongly complete list of sequences is when each $\sigma_i$ is $1$-complete for all $i$ where $1 \le i \le n$, i.e. each $\sigma_i$ contains all letters of $A$. For example, if $A=\{1,2,3\}$, then the list of sequences given by $\sigma_1=\langle 1, 2, 3 \rangle, \sigma_2=\langle 1, 2, 3 \rangle, \sigma_3=\langle 1, 2, 3 \rangle$ is both forward complete and backward complete, and therefore is strongly complete. In contrast, the list of sequences given by $\sigma_1=\langle 1, 2, 3 \rangle, \sigma_2=\langle 1, 2 \rangle, \sigma_3=\langle 1, 3 \rangle$ is forward complete but not backward complete, because $\langle 2 \rangle$ is not a subsequence of $\sigma_3$.

\begin{flushleft} \textbf{Definition 2.} A sequence $\sigma$ over $A$ is said to be a \textit{quasi-palindrome} if there exists a bijection $B: A \to A$ such that for all $k$ with $ 1 \le k \le l_{\sigma}$, $B(\sigma[k]) = \sigma[l_{\sigma}-k+1]$. Given a list of sequences $\sigma_1, \sigma_2,..., \sigma_n$, if $\sigma_1 \sigma_2 ... \sigma_n$ is a quasi-palindrome and for all $k$ with $ 1 \le k \le n$, $l_{\sigma_k} = l_{\sigma_{n-k+1}}$, then the list is said to be a quasi-palindrome.
\end{flushleft}

\begin{flushleft} In other words, a sequence is a quasi-palindrome if there is a bijection of $A$ to map a sequence to its reverse. For example $\langle 1,2,3,4  ,1,2,3, 4,1,2,  3,4,1,  2,3,4,1 \rangle$ is a quasi-palindrome as given by $B(1)=1$, $B(2)=4$, $B(3)=3$ and $B(4)=2$. Applying the bijection $B$ to individual elements of the sequence yields $\langle 1,4,3,2  ,1,4,3, 2,1,4,  3,2,1,  4,3,2,1 \rangle$, which is the reverse sequence of the original one. The following theorem illustrates the relationship between forward and backward completeness with quasi-palindrome.\end{flushleft}

\begin{flushleft} \textbf{Theorem 3.} 
Let $\sigma_1, \sigma_2, ..., \sigma_n$ be a list of sequences which is forward complete and is also a quasi-palindrome. Then the list is strongly complete.
\end{flushleft}
\begin{proof} 
Firstly, note that the reverse sequence of a $k$-complete sequence is still $k$-complete. This is because a $k$-complete sequence will contain every $k$-length sequence and its reverse as subsequence, so the reverse $k$-complete sequence will contain that too. To prove the theorem, we note that for any $k$, where $1 \le k \le n$, the bijection for the quasi-palindrome will map $\sigma_{n-k+1} \sigma_{n-k+2}... \sigma_n$ to the reverse sequence of $\sigma_1 \sigma_2 ... \sigma_k$. The later sequence is $k$-complete due to forward completeness, so the former sequence must be $k$-complete too due to the bijection. So the list is also backward complete and hence is strongly complete.
\end{proof}

The purpose of strong completeness is to enable construction of the supersequence over a larger set, based on the following theorem.

\begin{flushleft} \textbf{Theorem 4.} 
Let $\sigma_1, \sigma_2, ..., \sigma_n$ be a list of sequences over a set $A$ which is strongly complete, and $n = |A|$. If $x$ is a letter not in $A$, then  $x \sigma_1 x \sigma_2 x ... x \sigma_n x$ is a supersequence over $A \cup \{x\}$.
\end{flushleft}
\begin{proof} 
Given any permutation $\rho$ of all letters from $A \cup \{x\}$, if $\rho[1] = x$ or $\rho[-1] = x$, then $\rho < x \sigma_1 \sigma_2 ... \sigma_n$ or $\rho < \sigma_1 \sigma_2 ... \sigma_n x$ respectively, and hence $\rho < x \sigma_1 x \sigma_2 x ... x \sigma_n x$. Suppose $x$ is neither the first or the last element of $\rho$, so let $\rho = \langle a_1, a_2,... a_k, x, a_{k+1}, ... a_n\rangle$, where $1 \le k < n$. Then by definition of strongly complete, $\langle a_1, a_2,... a_k \rangle < \sigma_1 \sigma_2 ... \sigma_k$, and $\langle a_{k+1}, ... a_n \rangle < \sigma_{k+1} ... \sigma_n$. So $ \rho < \sigma_1 \sigma_2 ... \sigma_k x \sigma_{k+1...} \sigma_n < x \sigma_1 x \sigma_2 x ... x \sigma_n x$. The later is therefore a supersequence over $A \cup \{x\}$.
\end{proof}

Next, we will define the list of sequences at the first level called $T_1$. It can be used to construct the classical supersequences of Newey [7] and others. It is important also because it is the fundamental pattern that is part of higher level lists like $T_2, T_3$ and so on.

\begin{flushleft} \textbf{Definition 5.} Given an integer $n > 3$, we define $T_1(n)$ to be a list of $n$ many sequences $T_1(n,1), T_1(n,2),...,T_1(n,n)$ over $A=\{1,2,...,n\}$ as follows.
 \begin{enumerate}
  \item Define $T_1(n,1) = \sigma_1 = \langle1,2,...,n\rangle$ and $T_1(n,2) = \sigma_2 = \langle1,2,...,n-1\rangle$.
  \item For any $k$ where $2 < k < n$, define $T_1(n,k) = \sigma_k$ to be the sequence such that $l_{\sigma_k} = n-1$, $\sigma_k[1] = \sigma_{k-2}[-1]$, and for all $i$ where $ 2 \le i \le l_{\sigma_k}$, $\sigma_k[i] = \sigma_{k-1}[i-1]$. Equivalently, we have $$\sigma_k = \sigma_{k-2}[-1] \cdot \sigma_{k-1}[1,-2].$$
  \item Define $T_1(n,n) = \sigma_n$ to be the sequence such that $l_{\sigma_n} = n$, $\sigma_n[1] = \sigma_{n-2}[-1]$, and for all $i$ where $ 2 \le i \le l_{\sigma_n}$, $\sigma_n[i] = \sigma_{n-1}[i-1]$. Equivalently, we have $$\sigma_n = \sigma_{n-2}[-1] \cdot \sigma_{n-1}.$$
 \end{enumerate}
\end{flushleft}

Note that each of the $\sigma_k$ sequence in the $T_1$ list, except the first and last, omits a single letter from $A$.

\begin{flushleft} \textbf{Example 6.} The $T_1(6)$ list of sequences are defined as follows.
$T_1(6,1) = \sigma_1 = \langle1,2,3,4,5,6\rangle$.

$T_1(6,2) = \sigma_2 = \langle1,2,3,4,5\rangle$.

$T_1(6,3) = \sigma_3=\langle 6,1,2,3,4\rangle$.

$T_1(6,4) = \sigma_4=\langle 5,6,1,2,3\rangle$.

$T_1(6,5) = \sigma_5=\langle 4,5,6,1,2\rangle$.

$T_1(6,6) = \sigma_6 = \langle3,4,5,6,1,2\rangle$.

A bijection $B$ can be defined to demonstrate $T_1(6)$ is a quasi-palindrome as follows: $B(1)=2$, $B(2)=1$, $B(3)=6$, $B(4)=5$, $B(5)=4$ and $B(6)=3$. This is done by matching elements of $\sigma_1$ with elements of reverse of $\sigma_6$.
\end{flushleft}

The following two results will prove that $T_1$ is strongly complete. When combined with Theorem 4, this will provide a supersequence over $A \cup \{x\}$. The general strategy of the proof, which will later be expanded to cover higher level lists, is by induction on $k$ for $1 \le k \le n$. We will prove that for any $\rho \in [A]_k$, $\rho < \sigma_1 \sigma_2 ... \sigma_k$, thus establishing the $k$-completeness of the sequence and therefore the forward completeness of the list. The conclusion is obvious if $\sigma_k$ is $1$-complete because we can then use induction to get the result. So we will concentrate on cases when $\rho[k]$ is the missing letter in $\sigma_k$, i.e. $\rho[k] = \sigma_{k-1}[-1]$. This will entail recursively tracing backward to $\rho[k-1]$ or further, until we eventually get two consecutive elements of $\rho$ that come from the same $\sigma_i$ for some $i \le k$.

\begin{flushleft} \textbf{Lemma 7.} 
$T_1(n) = \sigma_1, \sigma_2, ..., \sigma_n$ is forward complete for each $n \ge 3$.
\end{flushleft}
\begin{proof} We prove that $\sigma_1 ... \sigma_k$ is $k$-complete for each $k \le n$ by induction on $k$. Clearly, $\sigma_1$ is $1$-complete. Suppose that $\sigma_1 ... \sigma_{k-1}$ is $(k-1)$-complete, and that $\rho$ is any sequence in $[A]_k$. If $\rho[k] \in \sigma_k$, then $\rho < \sigma_1 ... \sigma_k$ since $\rho|(k-1) < \sigma_1 ... \sigma_{k-1}$ by the induction hypothesis. On the other hand, if $\rho[k] \notin \sigma_k$, then $\rho[k] = \sigma_{k-1}[-1]$. Moreover, since $\rho|(k-1)$ doesn't include $\rho[k]$, we have $\rho|(k-1) < \sigma_1 ... \sigma_{k-2} \cdot \sigma_{k-1}[1,-2]$ by the induction hypothesis. Hence  $\rho < \sigma_1 ... \sigma_k$, which is thus $k$-complete as required.
\end{proof}

\begin{flushleft} \textbf{Theorem 8.} 
Given any integer $n > 3$, $T_1(n)$ is strongly complete over $A$.
\end{flushleft}
\begin{proof} 
The previous Lemma proves that $T_1(n)$ is forward complete. Backward completeness can be proven similarly but with induction going backward from $n$ down to $1$. Alternatively, it can be proven that the list of sequences is a quasi-palindrome. A bijection $B$ can be defined from the elements of $\sigma_1$ to the elements of the reverse of $\sigma_n$ as follows: for all i where $1 \le i \le n$, define $B(\sigma_1[i]) = \sigma_n[n-i+1]$. It can be verified that B will map $\sigma_2$ to the reverse of $\sigma_{n-1}$, $\sigma_3$ to the reverse of $\sigma_{n-2}$ and so on. Hence the list is a quasi-palindrome, and therefore strongly complete by Theorem 3.
\end{proof}

\section{Second Level list}

\begin{flushleft} Based on the previous presentation, we will define the list of sequences $T_2$ in this section. Using similar argument as before, we will then prove its strong completeness. When combined with Theorem 4, this list can be used to construct supersequence discovered by Radomirovic [8]. The proof provided here that expands from Theorem 8, however, is original. It is this new approach that enables generalization in later sections to yield additional saving with higher level lists. \end{flushleft}

\begin{flushleft} \textbf{Definition 9.} Given an integer $n$ with $n \ge 6$ and $n \equiv 0 \pmod{3}$, we define $T_2(n)$ to be a list of $n$ many sequences $T_2(n,1), T_2(n,2),...,T_2(n,n)$ over $A=\{1,2,...,n\}$ as follows.
 \begin{enumerate}
  \item Define $T_2(n,1) =T_1(n,1)$, $T_2(n,2) =T_1(n,2)$ and $T_2(n,3) =T_1(n,3)$.
  \item For any $k$ where $3 < k < n-2$ and $k \equiv 1 \pmod{3}$, define $T_2(n,k) = \sigma_k$ to be the sequence such that $l_{\sigma_k} = n-1$ as follows.
  \begin{enumerate}
    \item If $k=4$, then define  
    $$\sigma_k =\sigma_{k-2}[-1] \cdot \sigma_{k-1}[2,-3] \cdot n \cdot \sigma_{k-1}[-2] .$$
    \item If $k>4$, then define
    $$\sigma_k =\sigma_{k-2}[-1] \cdot \sigma_{k-1}[1] \cdot \sigma_{k-1}[3,-3] \cdot n \cdot \sigma_{k-1}[-2] .$$
   \end{enumerate}
   
  \item For any $k$ where $3 < k < n-2$ and $k \equiv 2 \pmod{3}$, define $T_2(n,k) = \sigma_k$ to be the sequence such that $l_{\sigma_k} = n-2$ and   
    $$\sigma_k =\sigma_{k-2}[-1] \cdot \sigma_{k-1}[1,-3].$$

  \item For any $k$ where $3 < k < n-2$ and $k \equiv 0 \pmod{3}$, define $T_2(n,k) = \sigma_k$ to be the sequence such that $l_{\sigma_k} = n-1$ and   
      $$\sigma_k =\sigma_{k-2}[-1] \cdot n \cdot \sigma_{k-1}[1,-2].$$

  \item Define 
    $$T_2(n,n-2) =\sigma_{n-2} = \sigma_{n-4}[-1] \cdot \sigma_{n-3}[1] \cdot \sigma_{n-3}[3,-2] \cdot n.$$
    $$T_2(n,n-1) =\sigma_{n-1} = \sigma_{n-3}[-1] \cdot \sigma_{n-2}[1,-2].$$
    $$T_2(n,n) =\sigma_{n} = \sigma_{n-2}[-1] \cdot \sigma_{n-1}= n \cdot \sigma_{n-1}.$$
 \end{enumerate}
\end{flushleft}

Note that each of the $\sigma_k$ sequence in the $T_2$ list defined above, except the first and last, omits one or two letters from $A$, with $\sigma_5, \sigma_8, \sigma_{11}, ... $ being the sequences that omit two letters.

\begin{flushleft} \textbf{Example 10.} The $T_2(12)$ list of sequences are given as follows. Compare with $T_1$, there is an additional skip letter $\textbf{12}$, which is bolded for ease of identification, removed in $\sigma_5$ and $\sigma_8$.

$\sigma_1 = \langle1,2,3,4,5,6,7,8,9,10,11,\textbf{12}\rangle$.

$\sigma_2 = \langle1,2,3,4,5,6,7,8,9,10,11\rangle$.

$\sigma_3 = \langle  \textbf{12},1,2,3,4,5,6,7,8,9,10\rangle$.

$\sigma_4=\langle 11,1,2,3,4,5,6,7,8, \textbf{12},9\rangle$.

$\sigma_5=\langle 10,11,1,2,3,4,5,6,7,8\rangle$.

$\sigma_6=\langle 9, \textbf{12},10,11,1,2,3,4,5,6,7\rangle$.

$\sigma_7=\langle 8,9,10,11,1,2,3,4,5, \textbf{12},6\rangle$.

$\sigma_8=\langle 7,8,9,10,11,1,2,3,4,5\rangle$.

$\sigma_9=\langle 6, \textbf{12},7,8,9,10,11,1,2,3,4\rangle$.

$\sigma_{10} = \langle5,6,7,8,9,10,11,1,2,3, \textbf{12}\rangle$.

$\sigma_{11} = \langle4,5,6,7,8,9,10,11,1,2,3\rangle$.

$\sigma_{12} = \langle \textbf{12},4,5,6,7,8,9,10,11,1,2,3\rangle$.

A bijection $B$ can be defined to demonstrate $T_2(12)$ is a quasi-palindrome as follows: $B(1)=3$, $B(2)=2$, $B(3)=1$, $B(4)=11$, $B(5)=10$, $B(6)=9$, $B(7)=8$, $B(8)=7$, $B(9)=6$, $B(10)=5$, $B(11)=4$ and $B(12)=12$. This is done again by matching elements of $\sigma_1$ with elements of reverse of $\sigma_{12}$. It can also be verified that $B(\sigma_2)$ is the reverse of $\sigma_{11}$, $B(\sigma_3)$ is the reverse of $\sigma_{10}$, and so on.
\end{flushleft}

The key improvement of $T_2$ over $T_1$ is to skip additional letter $n$ at those sequences $\sigma_k$ where $k \equiv 2 \pmod{3}$. We introduce additional notation and a key Lemma, and then prove the strong completeness of the newly defined list. The proof needs to show that for any $k$-length sequence $\rho$, $\rho < \sigma_1 \sigma_2 ... \sigma_k$. As before, this is difficult only if $\rho[k] \notin \sigma_k$. In this case, we again trace the elements of $\rho$ recursively backward within each $\sigma_i$ for $i \le k$. We succeed when we find two consecutive elements of $\rho$ that come from the same $\sigma_i$.

\begin{flushleft} \textbf{Definition 11.} Let $\sigma$ be any sequence of distinct letters. For any $a \in \sigma$, let $\sigma[>a]$ represents the set of elements in $\sigma$ that comes later than $a$. Formally, we define $$\sigma[>a] = \{b \in \sigma: \sigma^{-1}[b] > \sigma^{-1}[a] \}.$$ 
\end{flushleft}

\begin{flushleft} \textbf{Lemma 12.} Suppose $\sigma_1, \sigma_2, ..., \sigma_n$ is a $T_2(n)$ list of sequences over $A$. Let $k$ be any integer such that $3 < k < n-2$ and $k \equiv 2 \pmod{3}$. Assume that  $\sigma_1, \sigma_2, ... ,\sigma_{k-2}$ is forward complete. If $\rho \in [A]_{k-1}$, $n \notin \rho$ and $\rho[k-1]=\sigma_{k-1}[-1]$, then $\rho < \sigma_1 \sigma_2 ... \sigma_{k-2}$.
\end{flushleft}
\begin{proof} The proof follows the style of Lemma 7, with additional assumption that $n \notin \rho$. In the following, note that $\rho[k-1] = \sigma_{k-1}[-1] = \sigma_{k-2}[-2]$. Consider the various possible cases of the next element $\rho[k-2]$.

\begin{enumerate}
  \item $\rho[k-2] = \sigma_{k-3}[-1]$. Apply the argument of Lemma 7 to $\rho|k-2$, which shows that it is a subsequence of $\sigma_1 \sigma_2 ... \sigma_{k-3}$. Since $\rho[k-1] \in \sigma_{k-2}$, we have $\rho < \sigma_1 \sigma_2 ... \sigma_{k-2}$
  \item $\rho[k-2] \in \sigma_{k-2}$ and $\rho[k-2] \notin \sigma_{k-2}[>\rho[k-1]] = \{\sigma_{k-2}[-1] \}$. Then $\langle \rho[k-2], \rho[k-1] \rangle < \sigma_{k-2}$. Together with the fact that $\sigma_1, \sigma_2, ... ,\sigma_{k-3}$ is $(k-3)$-complete, the Lemma is true.
  \item $\rho[k-2] \in \sigma_{k-2}[>\rho[k-1]] = \{\sigma_{k-2}[-1] \}$, i.e. $\rho[k-2]$ is the last element of $\sigma_{k-2}$, then recursively consider the next element $\rho[k-3]$.
\end{enumerate}

If the process continues to $\sigma_1$, then $\rho$ is having the last element of each $\sigma_i$ in the order. However, the last element of $\sigma_1$ is $n$, which is not an element of $\rho$ by assumption. So, $\langle \rho[1], \rho[2] \rangle < \sigma_1$. This implies that the Lemma is true.
\end{proof}

Fundamentally, the reason that the Lemma is true is because there are only $k-2$ integers, excluding $n$, greater or equal to $n-k+2 = \sigma_{k-1}[-1]$. But since $\rho$ is of length $k-1$, at least a member of $\rho$ is not the last element of some $\sigma_i$. So, there will be at least two consecutive elements of $\rho$ that come from the same $\sigma_i$.

\begin{flushleft} \textbf{Theorem 13.} 
Given an integer $n \ge 6$ such that $n \equiv 0 \pmod{3}$, $T_2(n)$ is strongly complete over $A$.
\end{flushleft}
\begin{proof} 
For most of $k \le n$, the proof for $T_2(n,k)$ is exactly the same as Theorem 8. The only additional case to consider is for $k$ where $3 < k < n-2$ and $k \equiv 2 \pmod{3}$. This is when the sequence skips additional letter $n$, i.e. case 3 in Definition 9 and Example 10. In that case, note that the only two letters from $A$ that are missing from the elements of $\sigma_k$ are $n$ and $n-k+2$. The later is the last element of $\sigma_{k-1}$, therefore, any sequence $\rho \in [A]_k$ with $\rho[k] = n-k+2$ is a subsequence of $\sigma_1 \sigma_2 ... \sigma_{k-1}$ as in the proof of Lemma 7. So it remains to be proven that $\forall \rho \in [A]_k$ with $\rho[k] = n$, it is true that $\rho < \sigma_1 \sigma_2 ... \sigma_k$. Given such $\rho$, the following will prove a stronger result that $\rho < \sigma_1 \sigma_2 ... \sigma_{k-1}$.

Let $\sigma' = \sigma_{k-1} | n-3$ denote the initial substring of $\sigma_{k-1}$ except the last two elements. Recall the last two elements are exactly $n$ and $n-k+2$, the letters that are missing in $\sigma_k$. We first prove the following Claim.

\begin{flushleft} \textbf{Claim 14.} 
The sequence $\rho | k-1$ is a subsequence of $\sigma_1 \sigma_2 ... \sigma_{k-2} \sigma'$.
\end{flushleft}
\begin{proof} 
By induction, $\sigma_1 \sigma_2 ... \sigma_{k-2}$ is $(k-2)$-complete, so $\rho | k-2 < \sigma_1 \sigma_2 ... \sigma_{k-2}$. If $\rho[k-1]$ is an element of $\sigma'$, then the Claim is true. There are only three letters of $A$ which are not an element of $\sigma'$. They are $n$, $n-k+2$, and $n-k+3$. So we only need to consider if $\rho[k-1]$ is one of them.
\begin{enumerate}
  \item $\rho[k-1] = n$. This is not possible because $\rho[k] = n$.
  \item $\rho[k-1] = n-k+2 = \sigma_{k-1}[-1]$. Apply Lemma 12 to $\rho|k-1$ to prove that it is a subsequence of $\sigma_1 \sigma_2 ... \sigma_{k-2}$.
  \item $\rho[k-1] = n-k+3 = \sigma_{k-2}[-1]$. Then $\rho | k-1 < \sigma_1 \sigma_2 ... \sigma_{k-2}$ as in the proof of Lemma 7.
\end{enumerate}
\end{proof}

Since $\rho[k] = n$ and $\sigma_{k-1}[-2] = n$, so $\rho < \sigma_1 \sigma_2 ... \sigma_{k-1}$. Therefore $ \sigma_1 \sigma_2 ... \sigma_k$ is $k$-complete. The backward completeness can either be proven with the above argument in reverse or using the quasi-palindrome property of $T_2(n)$.
\end{proof}

\begin{flushleft} \textbf{Example 15.} Using the sequences defined in Example 10, we will show that $\sigma_1 \sigma_2 ... \sigma_8$ is $8$-complete. Since $\sigma_1 \sigma_2 ... \sigma_7$ is $7$-complete, any $\rho \in [A]_8$ with $\rho[8] \ne 12$ or $\rho[8] \ne 6$ (the only two missing letters in $\sigma_8$) can be easily seen to be a subsequence of $\sigma_1 \sigma_2 ... \sigma_8$. For $\rho[8] = 6$, then $\rho < \sigma_1 \sigma_2 ... \sigma_7$ by Lemma 7, because $6$ is the last element of $\sigma_7$. So, what remains is to show that if $\rho[8]=12$, than it is also a subsequence of $\sigma_1 \sigma_2 ... \sigma_8$. Claim 14 asserts that such $\rho|7$ can be generated before the last two elements of $\sigma_7$. If $\rho[7] \ne 12$, $\rho[7] \ne 6$ or $\rho[7] \ne 7$, then the Claim is true because $\rho[7]$ will be an element in $\sigma_7$ before the last two elements. Since $\rho[8] = 12$, so $\rho[7] \ne 12$. Lemma 12 and Lemma 7 shows that if $\rho[7] = 6$ or $\rho[7] = 7$, then $\rho |7 < \sigma_1 \sigma_2 ... \sigma_6$. Therefore $\rho < \sigma_1 \sigma_2 ... \sigma_7$ and $\sigma_1 \sigma_2 ... \sigma_8$ must be $8$-complete.
\end{flushleft}

\section{Higher Level lists}
\begin{flushleft} This section generalizes previous definitions to define sequences with additional skip letters. We assume $s$ to be an integer $\ge 3$ for this section and use it to index the different levels of list. For each $T_s$ to be defined, there will be a total of $s$ many elements removed in the appropriate sequence $\sigma_k$. However, only $s-1$ of them are called the skip letters as defined in Definition 16, because they are the same letters that can be removed in different cycles. The last element of the previous sequence $\sigma_{k-1}$, which will always be removed like in all places, will not be called the skip letter.

\begin{flushleft} \textbf{Definition 16.} Given integer $n$, define the $s-1$ many integers $n-s+2,...,n-1,n$ to be the \textit{skip letters} for the $T_s(n)$ list of sequences. We use $\phi_s$ and $\phi_{-s}$ to denote the sequence $\langle n-s+2,...,n-1,n \rangle$ and its reverse sequence $\langle n, n-1,..., n-s+2 \rangle$ respectively.
\end{flushleft} 

Definition 17 defines $T_s$. It is interleaved with an example to help understanding of the sequence pattern. The cases of Definition 17 can be divided into different categories: cases 1 and 7, each with $s+1$ sequences, are the initial and final stages; cases 2 and 5 are where skip letters jump to the end or start of the sequences respectively; cases 3 and 6 are normal $T_1$ style forwarding; and finally case 4 is where all the skip letters are removed that provides the saving. A complete cycle consists of sequences from case 2, 3, 4, 5 and 6 with a total of $2s-1$ sequences.
\end{flushleft}

\begin{flushleft} \textbf{Definition 17.} Given an integer $n$ with $n \ge 4s+1$ and $n \equiv 3 \pmod{2s-1}$, we define $T_s(n)$ to be a list of $n$ many sequences $T_s(n,1), T_s(n,2),...,T_s(n,n)$ over $A=\{1,2,...,n\}$ as follows.
 \begin{enumerate}
  \item For integer $k$, $1\le k \le s+1$, we define $T_s(n,k) = \sigma_k = T_1(n,k)$. 
  
  \begin{flushleft} \textbf{Example 18.1.} Assume $s=3, n=18$ and define $T_3(18)$ as follows. The two skip letters $\textbf{17}$ and $\textbf{18}$ are bolded for ease of reading. The purpose of these $4$ sequences is to wait for all skip letters to go to the beginning of the sequences.
 $T_3(18,1) = \sigma_1 = \langle1,2,3,4,5,6,7,8,9,10,11,12,13,14,15,16,\textbf{17},\textbf{18}\rangle$,
 $T_3(18,2) = \sigma_2 = \langle1,2,3,4,5,6,7,8,9,10,11,12,13,14,15,16,\textbf{17}\rangle$,
 $T_3(18,3) = \sigma_3 = \langle \textbf{18},1,2,3,4,5,6,7,8,9,10,11,12,13,14,15,16\rangle$,
 $T_3(18,4) = \sigma_4 = \langle \textbf{17},\textbf{18},1,2,3,4,5,6,7,8,9,10,11,12,13,14,15\rangle$.
\end{flushleft}

  \item For any $k$ where $s+1 < k < n-s$ and $k \equiv s+2 \pmod{2s-1}$, define $T_s(n,k) = \sigma_k$ to be a sequence such that $l_{\sigma_k} = n-1$ with the following subcases.
  \begin{enumerate}
    \item If $k=s+2$, then
    $$\sigma_k =\sigma_{k-2}[-1] \cdot \sigma_{k-1}[s,-s-1] \cdot \phi_s \cdot \sigma_{k-1}[-s,-2].$$
    \item If $k>s+2$, then
    $$\sigma_k =\sigma_{k-2}[-1] \cdot \sigma_{k-1}[1,s-1] \cdot \sigma_{k-1}[2s-1,-s-1] \cdot \phi_s \cdot \sigma_{k-1}[-s,-2] .$$
   \end{enumerate}
   
  \begin{flushleft} \textbf{Example 18.2.} All skip letters are jumped to the last $s$-th position.
 $T_3(18,5) =\sigma_5 = \langle16,1,2,3,4,5,6,7,8,9,10,11,12,\textbf{17},\textbf{18},13,14\rangle$
 $T_3(18,10) =\sigma_{10} = \langle11,12,13,14,15,16,1,2,3,4,5,6,7,\textbf{17},\textbf{18},8,9\rangle$
\end{flushleft}
   
  \item For any $k$ where $s+1 < k < n-s$. Let $k$ range from $k \equiv s+3 \pmod{2s-1}$ to $k \equiv 1 \pmod{2s-1}$. Define $T_s(n,k) = \sigma_k$ to be a sequence such that $l_{\sigma_k} = n-1$ and
    $$\sigma_k =\sigma_{k-2}[-1] \cdot \sigma_{k-1}[1,-2].$$

  \begin{flushleft} \textbf{Example 18.3.} $T_1$ style forwarding till skip letters are ready at the end.
 $T_3(18,6) =\sigma_6 = \langle15,16,1,2,3,4,5,6,7,8,9,10,11,12, \textbf{17},\textbf{18},13\rangle$
 $T_3(18,11) =\sigma_{11} = \langle10,11,12,13,14,15,16,1,2,3,4,5,6,7,\textbf{17},\textbf{18},8\rangle$
\end{flushleft}

  \item For any $k$ where $s+1 < k < n-s$ and $k \equiv 2 \pmod{2s-1}$. Define $T_s(n,k) = \sigma_k$ to be a sequence such that $l_{\sigma_k} = n-1$ and
    $$\sigma_k =\sigma_{k-2}[-1] \cdot \sigma_{k-1}[1,-s-1].$$
    
  \begin{flushleft} \textbf{Example 18.4.} All skip letters are removed.
 $T_3(18,7) =\sigma_7 = \langle14,15,16,1,2,3,4,5,6,7,8,9,10,11,12\rangle$
 $T_3(18,12) =\sigma_{12} = \langle9,10,11,12,13,14,15,16,1,2,3,4,5,6,7\rangle$
\end{flushleft}

  \item For any $k$ where $s+1 < k < n-s$ and $k \equiv 3 \pmod{2s-1}$. Define $T_s(n,k) = \sigma_k$ to be a sequence such that $l_{\sigma_k} = n-1$ and   
    $$\sigma_k =\sigma_{k-2}[-1] \cdot \phi_{-s} \cdot \sigma_{k-1}[1,-2].$$
    
  \begin{flushleft} \textbf{Example 18.5.} Skip letters are recovered in reverse order and placed in the start of sequences to ensure quasi-palindrome property.
 $T_3(18,8) =\sigma_8 = \langle13,\textbf{18},\textbf{17},14,15,16,1,2,3,4,5,6,7,8,9,10,11\rangle$
 $T_3(18,13) =\sigma_{13} = \langle8,\textbf{18},\textbf{17},9,10,11,12,13,14,15,16,1,2,3,4,5,6\rangle$
\end{flushleft}

  \item For any $k$ where $s+1 < k < n-s$. Let $k$ range from $k \equiv 4 \pmod{2s-1}$ to $k \equiv s+1 \pmod{2s-1}$. Define $T_s(n,k) = \sigma_k$ to be a sequence such that $l_{\sigma_k} = n-1$ and   
    $$\sigma_k =\sigma_{k-2}[-1] \cdot \sigma_{k-1}[1,-2].$$

  \begin{flushleft} \textbf{Example 18.6.} $T_1$ style forwarding for the next skip cycle (if available).
 $T_3(18,9) =\sigma_9 = \langle12,13,\textbf{18},\textbf{17},14,15,16,1,2,3,4,5,6,7,8,9,10\rangle$
 $T_3(18,14) =\sigma_{14} = \langle7,8,\textbf{18},\textbf{17},9,10,11,12,13,14,15,16,1,2,3,4,5\rangle$
\end{flushleft}

  \item Define 
    $$T_s(n,n-s) =\sigma_{n-s} = \sigma_{n-s-2}[-1] \cdot \sigma_{n-s-1}[1,s-1] \cdot \sigma_{n-s-1}[2s-1,-2] \cdot \phi_{-s}.$$
    For any $i$ where $1 \le i < s$, define
    $$T_s(n,n-s+i) =\sigma_{n-s+i} = \sigma_{n-s+i-2}[-1] \cdot \sigma_{n-s+i-1}[1,-2].$$
    And finally, define
    $$T_s(n,n) =\sigma_{n} = \sigma_{n-2}[-1] \cdot \sigma_{n-1}.$$
    
  \begin{flushleft} \textbf{Example 18.7.} These $4$ sequences is to ensure quasi-palindrome by reversing the initial pattern of the $4$ sequences in Example 18.1.
 $T_3(18,15) =\sigma_{15} = \langle6,7,8,9,10,11,12,13,14,15,16,1,2,3,4,\textbf{18},\textbf{17}\rangle$
 $T_3(18,16) =\sigma_{16} = \langle5,6,7,8,9,10,11,12,13,14,15,16,1,2,3,4,\textbf{18}\rangle$
 $T_3(18,17) =\sigma_{17} = \langle \textbf{17},5,6,7,8,9,10,11,12,13,14,15,16,1,2,3,4\rangle$
 $T_3(18,18) =\sigma_{18} = \langle \textbf{18},\textbf{17},5,6,7,8,9,10,11,12,13,14,15,16,1,2,3,4\rangle$
\end{flushleft}

 \end{enumerate}
\end{flushleft}
Note that each of the $\sigma_k$ sequence in the $T_s$ list defined above, except the first and last, omits either one letter or $s$ letters from $A$, with $\sigma_{2s+1}, \sigma_{4s}, \sigma_{6s-1}, ... $ being the sequences that omit $s$ letters.

\begin{flushleft} \textbf{Example 18.} The whole example of $T_3(18)$ is being reproduced below. 
 $\sigma_1 = \langle1,2,3,4,5,6,7,8,9,10,11,12,13,14,15,16,\textbf{17},\textbf{18}\rangle$.
 $\sigma_2 = \langle1,2,3,4,5,6,7,8,9,10,11,12,13,14,15,16,\textbf{17}\rangle$.
 $\sigma_3 = \langle \textbf{18},1,2,3,4,5,6,7,8,9,10,11,12,13,14,15,16\rangle$.
 $\sigma_4 = \langle \textbf{17},\textbf{18},1,2,3,4,5,6,7,8,9,10,11,12,13,14,15\rangle$.
 $\sigma_5 = \langle16,1,2,3,4,5,6,7,8,9,10,11,12,\textbf{17},\textbf{18},13,14\rangle$.
 $\sigma_6 = \langle15,16,1,2,3,4,5,6,7,8,9,10,11,12, \textbf{17},\textbf{18},13\rangle$.
 $\sigma_7 = \langle14,15,16,1,2,3,4,5,6,7,8,9,10,11,12\rangle$.
 $\sigma_8 = \langle13,\textbf{18},\textbf{17},14,15,16,1,2,3,4,5,6,7,8,9,10,11\rangle$.
 $\sigma_9 = \langle12,13,\textbf{18},\textbf{17},14,15,16,1,2,3,4,5,6,7,8,9,10\rangle$.
 $\sigma_{10} = \langle11,12,13,14,15,16,1,2,3,4,5,6,7,\textbf{17},\textbf{18},8,9\rangle$.
 $\sigma_{11} = \langle10,11,12,13,14,15,16,1,2,3,4,5,6,7,\textbf{17},\textbf{18},8\rangle$.
 $\sigma_{12} = \langle9,10,11,12,13,14,15,16,1,2,3,4,5,6,7\rangle$.
 $\sigma_{13} = \langle8,\textbf{18},\textbf{17},9,10,11,12,13,14,15,16,1,2,3,4,5,6\rangle$.
 $\sigma_{14} = \langle7,8,\textbf{18},\textbf{17},9,10,11,12,13,14,15,16,1,2,3,4,5\rangle$.
 $\sigma_{15} = \langle6,7,8,9,10,11,12,13,14,15,16,1,2,3,4,\textbf{18},\textbf{17}\rangle$.
 $\sigma_{16} = \langle5,6,7,8,9,10,11,12,13,14,15,16,1,2,3,4,\textbf{18}\rangle$.
 $\sigma_{17} = \langle \textbf{17},5,6,7,8,9,10,11,12,13,14,15,16,1,2,3,4\rangle$.
 $\sigma_{18} = \langle \textbf{18},\textbf{17},5,6,7,8,9,10,11,12,13,14,15,16,1,2,3,4\rangle$.

The bijection that maps elements of $\sigma_1$ to the elements of the reverse of $\sigma_{18}$ will demonstrate the quasi-palindrome of the sequences: $B(1)=4$, $B(2)=3$, $B(3)=2$, $B(4)=1$, $B(5)=16$, $B(6)=15$, $B(7)=14$, $B(8)=13$, $B(9)=12$, $B(10)=11$, $B(11)=10$, $B(12)=9$, $B(13)=8$, $B(14)=7$, $B(15)=6$, $B(16)=5$, $B(17)=17$ and $B(18)=18$,
\end{flushleft}

The following Lemma contains the core technical proof of the paper.

\begin{flushleft} \textbf{Lemma 19.} Suppose $\sigma_1, \sigma_2, ..., \sigma_n$ is a $T_s(n)$ list of sequences over $A$. Given integer $k$ where $s+1 < k < n-s$ and $k \equiv 2 \pmod{2s-1}$. Assume $\rho \in [A]_k$ and $\rho[k]$ is one of the skip letters. Suppose $ \sigma_1 \sigma_2 ... \sigma_{k-1}$ is forward complete. Then $\rho < \sigma_1 \sigma_2 ... \sigma_{k-1}$.
\end{flushleft}
\begin{proof} The proof is similar to Lemma 12, but instead of allowing only the last element of a sequence to be the only possible entry to the next round of recursion, this proof will allow a set $M_i$ of size no bigger than $s-1$ many integers to be the entries. We will examine the possible occurrence of $\rho[i], i \le k-1$ in the $\sigma_j$ for $j \le i$. The following demonstrates the recursion up to two rounds.
  \begin{enumerate}
    \item $\rho[k-1] = \sigma_{k-2}[-1]$, then by Lemma 7 argument, $\rho|k-1 < \sigma_1 \sigma_2 ... \sigma_{k-2}$, and since $\rho[k] \in \sigma_{k-1}$, therefore $\rho < \sigma_1 \sigma_2 ... \sigma_{k-1}$.
    \item $\rho[k-1] \in \sigma_{k-1}$ and $\rho[k-1] \notin M_{k-1}$, where $M_{k-1} = \sigma_{k-1}[>\rho[k]]$. then $\langle \rho[k-1], \rho[k] \rangle < \sigma_{k-1}$, Together with the $(k-2)$-completeness of $ \sigma_1 \sigma_2 ... \sigma_{k-2}$, these imply that $\rho < \sigma_1 \sigma_2 ... \sigma_{k-1}$.
    
    \item $\rho[k-1] \in M_{k-1}$ and continue the process by considering $\rho[k-2]$.
      \begin{enumerate}
        \item $\rho[k-2] = \sigma_{k-3}[-1]$. The same argument applies.
        \item $\rho[k-2] \in \sigma_{k-2}$ and $\rho[k-2] \notin M_{k-2}$, where $M_{k-2} = \sigma_{k-2}[>\rho[k-1]] \setminus \{ \rho[k]\}$. So $\langle \rho[k-2], \rho[k-1] \rangle < \sigma_{k-2}$, Together with the $(k-3)$-completeness of $ \sigma_1 \sigma_2 ... \sigma_{k-3}$, these imply that $\rho < \sigma_1 \sigma_2 ... \sigma_{k-1}$.
         \item $\rho[k-2] \in M_{k-2}$ and the recursion continues.
     \end{enumerate}
    \end{enumerate}

At each step $k-i$, we remove all elements that previously occur in $\rho$, i.e. $M_{k-i} = \sigma_{k-i}[> \rho[k-i+1]] \setminus \{ \rho[k], \rho[k-1], ..., \rho[k-i+2]\}$. This has the effect of removing only the skip letters that occur in $\rho$. It can be shown that by the time $i$ reaches the case such that $k-i \equiv s+2 \pmod{2s-1}$ (i.e. case 2 of Definition 17, where all the skip letters are at the last $s$-th positon of the sequence), none of the skip letter is in $M_{k-i}$ (though the skip letters may appear in $M_{k-i}$ again for future cycle). That is because all the skip letters will either have occurred in $\rho$ and therefore explicitly removed from the definition of $M_{k-i}$, or some $\rho[k-j]$ with $j < i$ will have position $>$ the positons of all skip letters in $\sigma_{k-j}$. This is critical because it ensures that no elements with lower position than the skip letters in $\sigma_{k-i}$ can get into $M_{k-i-1}$, thus limiting the number of integers in it. The number of elements of $M_{k-i}$ plus the number of skip letters that have not occurred in $\rho$ previously is at most $s-1$. So the number of elements of $M_{k-i}$ is always $\le s-1$.

The number of integers in $M_{k-i}$ will be reduced to zero at or before $\sigma_1$. When that happens, the recursion terminates. If the recursion continues to $\sigma_1$, it is always true that $\langle \rho[1], \rho[2] \rangle < \sigma_1$. So the Lemma is proven.
\end{proof}

\begin{flushleft} \textbf{Example 20.} Consider $k=12$ for Example 18. If $\rho \in [A]_{12}$ and $\rho[12] = 17$. The recursion is as follows.
  \begin{enumerate}
    \item $\rho[11] = \sigma_{10}[-1] = 9$, then by Lemma 7, $\rho|11 < \sigma_1 \sigma_2 ... \sigma_{10}$, and since $\rho[12]=17 \in \sigma_{11}$, therefore $\rho < \sigma_1 \sigma_2 ... \sigma_{11}$.
    \item $\rho[11] \in \sigma_{11}$ and $\rho[11] \notin M_{11}$, where $M_{11} = \sigma_{11}[>17] = \{18,8\}$. then $\langle \rho[11], 17 \rangle < \sigma_{11}$, Together with the $10$-completeness of $ \sigma_1 \sigma_2 ... \sigma_{10}$, these imply that $\rho < \sigma_1 \sigma_2 ... \sigma_{11}$.
    
    \item $\rho[11] \in M_{11} = \{18,8\}$ and continue the process by considering $\rho[10]$.
      \begin{enumerate}
        \item $\rho[10] = \sigma_{9}[-1] = 10$. The same argument applies.
        \item $\rho[10] \in \sigma_{10}$ and $\rho[10] \notin M_{10}$, where $M_{10} = \sigma_{10}[>\rho[11]] \setminus \{17\}$.So $\langle \rho[10], \rho[11] \rangle < \sigma_{10}$, Together with the $9$-completeness of $ \sigma_1 \sigma_2 ... \sigma_9$, these imply that $\rho < \sigma_1 \sigma_2 ... \sigma_{11}$. 
         \item $\rho[10] \in M_{10}$, and the recursion continues.
     \end{enumerate}
    \end{enumerate}

Note that if $\rho[11] = 18$, then $M_{10} = \{8,9\}$, and if $\rho[11] = 8$, then $M_{10} = \{9\}$. In either case, no skip letter is in $M_{10}$. This is important because otherwise, a lot of elements will be qualified to be in $M_9$. Assume the maximum number of elements for each $M_i$, the following will be the sets $M_{10} = \{8,9\}$, $M_9 = \{9,10\}$, $M_8 = \{10,11\}$, $M_7 = \{11,12\}$, $M_6 = \{12,13\}$ (note if $18$ has not occurred previously, i.e.  $\rho[11] \ne 18$, then $M_6 = \{18,13\}$ instead), $M_5 = \{13,14\}$, $M_4 = \{14,15\}$, $M_3 = \{15,16\}$, $M_2 = \{16\}$, $M_1 = \{\}$. Note that the union of all $M_i$ sets is just the set of integers greater or equal to $8$, excluding $17 = \rho[12]$.
\end{flushleft}

Establishing the strong completeness of $T_s(n)$ concludes the necessary generalization.

\begin{flushleft} \textbf{Theorem 21.} $T_s(n)$ is strongly complete over $A$.
\end{flushleft}
\begin{proof} As in Theorem 13, the only case that needs to be considered is for sequences that have skip letters removed, i.e. case 4 in Definition 17. Let $k$ be $s+1 < k < n-s$ and $k \equiv 2 \pmod{2s-1}$. Let $\rho \in [A]_k$, and we need to prove that $\rho < \sigma_1 \sigma_2 ... \sigma_k$. The argument is the same as previously if $\rho[k]$ is not one of the skip letters. So assume $\rho[k]$ is a skip letter and using Lemma 19 to conclude that $\rho < \sigma_1 \sigma_2... \sigma_{k-1}$. Therefore $T_s(n)$ is forward complete. The proof for backward completenss is the same or can be proven using quasi-palindrome property of $T_s(n)$.
\end{proof}

\begin{flushleft} \textbf{Example 22.} The list $T_4(24)$ with $s=4, n=24$ is provided below.

 $\sigma_1 = \langle1,2,3,4,5,6,7,8,9,10,11,12,13,14,15,16,17,18,19,20,21,\textbf{22},\textbf{23},\textbf{24}\rangle$,
 $\sigma_2 = \langle1,2,3,4,5,6,7,8,9,10,11,12,13,14,15,16,17,18,19,20,21,\textbf{22},\textbf{23}\rangle$,
 $\sigma_3 = \langle\textbf{24},1,2,3,4,5,6,7,8,9,10,11,12,13,14,15,16,17,18,19,20,21,\textbf{22}\rangle$,
 $\sigma_4 = \langle\textbf{23},\textbf{24},1,2,3,4,5,6,7,8,9,10,11,12,13,14,15,16,17,18,19,20,21\rangle$,
 $\sigma_5 = \langle\textbf{22},\textbf{23},\textbf{24},1,2,3,4,5,6,7,8,9,10,11,12,13,14,15,16,17,18,19,20\rangle$,
 $\sigma_6 = \langle21,1,2,3,4,5,6,7,8,9,10,11,12,13,14,15,16,\textbf{22},\textbf{23},\textbf{24},17,18,19\rangle$,
 $\sigma_7 = \langle20,21,1,2,3,4,5,6,7,8,9,10,11,12,13,14,15,16,\textbf{22},\textbf{23},\textbf{24},17,18\rangle$,
 $\sigma_8 = \langle19,20,21,1,2,3,4,5,6,7,8,9,10,11,12,13,14,15,16,\textbf{22},\textbf{23},\textbf{24},17\rangle$,
 $\sigma_9 = \langle18,19,20,21,1,2,3,4,5,6,7,8,9,10,11,12,13,14,15,16\rangle$,
 $\sigma_{10} = \langle17,\textbf{24},\textbf{23},\textbf{22},18,19,20,21,1,2,3,4,5,6,7,8,9,10,11,12,13,14,15\rangle$,
 $\sigma_{11} = \langle16,17,\textbf{24},\textbf{23},\textbf{22},18,19,20,21,1,2,3,4,5,6,7,8,9,10,11,12,13,14\rangle$,
 $\sigma_{12} = \langle15,16,17,\textbf{24},\textbf{23},\textbf{22},18,19,20,21,1,2,3,4,5,6,7,8,9,10,11,12,13\rangle$,
 $\sigma_{13} = \langle14,15,16,17,18,19,20,21,1,2,3,4,5,6,7,8,9,\textbf{22},\textbf{23},\textbf{24},10,11,12\rangle$,
 $\sigma_{14} = \langle13,14,15,16,17,18,19,20,21,1,2,3,4,5,6,7,8,9,\textbf{22},\textbf{23},\textbf{24},10,11\rangle$,
 $\sigma_{15} = \langle12,13,14,15,16,17,18,19,20,21,1,2,3,4,5,6,7,8,9,\textbf{22},\textbf{23},\textbf{24},10\rangle$,
 $\sigma_{16} = \langle11,12,13,14,15,16,17,18,19,20,21,1,2,3,4,5,6,7,8,9\rangle$,
 $\sigma_{17} = \langle10,\textbf{24},\textbf{23},\textbf{22},11,12,13,14,15,16,17,18,19,20,21,1,2,3,4,5,6,7,8\rangle$,
 $\sigma_{18} = \langle9,10,\textbf{24},\textbf{23},\textbf{22},11,12,13,14,15,16,17,18,19,20,21,1,2,3,4,5,6,7\rangle$,
 $\sigma_{19} = \langle8,9,10,\textbf{24},\textbf{23},\textbf{22},11,12,13,14,15,16,17,18,19,20,21,1,2,3,4,5,6\rangle$,
 $\sigma_{20} = \langle7,8,9,10,11,12,13,14,15,16,17,18,19,20,21,1,2,3,4,5,\textbf{24},\textbf{23},\textbf{22}\rangle$,
 $\sigma_{21} = \langle6,7,8,9,10,11,12,13,14,15,16,17,18,19,20,21,1,2,3,4,5,\textbf{24},\textbf{23}\rangle$,
 $\sigma_{22} = \langle\textbf{22},6,7,8,9,10,11,12,13,14,15,16,17,18,19,20,21,1,2,3,4,5,\textbf{24}\rangle$,
 $\sigma_{23} = \langle\textbf{23},\textbf{22},6,7,8,9,10,11,12,13,14,15,16,17,18,19,20,21,1,2,3,4,5\rangle$,
 $\sigma_{24} = \langle\textbf{24},\textbf{23},\textbf{22},6,7,8,9,10,11,12,13,14,15,16,17,18,19,20,21,1,2,3,4,5\rangle$.\end{flushleft}
 
 There are a total of $548$ elements in all the $24$ sequences. Using Theorem 4 to interpose $25$ occurrences of a new element say $x$, we obtain the supersequence  $x \sigma_1 x \sigma_2 x ... x \sigma_{24} x$ over a set of $25$ letters. The size of the supersequence is $548+25=573$.

\section{Size of Supersequence}

\begin{flushleft} With the main tool proven, we can formally calculate the number of elements of the supersequence constructed through $T_s(n)$.
\end{flushleft}

\begin{flushleft} \textbf{Theorem 23.} Given integer $n \ge 4s+1$ and $n \equiv 3 \pmod{2s-1}$, there exists a supersequence of size $m^2 - \frac {5s-3}{2s-1}m +\frac {2s^2+9s-7}{2s-1}$ over a set of $m=n+1$ letters.
\end{flushleft}
\begin{proof} Using Theorem 4, a supersequence over a set of $m=n+1$ letters can be constructed from $T_s(n)$. This is done by interposing $m$ many occurrences of the new letter to $\sigma_1 \sigma_2... \sigma_n$. So the length of the supersequence is  $m + l_{\sigma_1 \sigma_2... \sigma_n}$.

Let $t = (m-2s-3)/(2s-1)$ be the number of skip sequences $\sigma_k$ where $k \equiv 2 \pmod{2s-1}$, i.e. the number of sequences defined in case 4 of Definition 17. Then $t$ also represents the number of cycles of $(2s-1)$-many sequences in the middle of $T_s(n)$. To calculate the number of elements of $\sigma_1 \sigma_2... \sigma_n$, the length of each $\sigma_i$ is added to get the following:
$$ l_{\sigma_1 \sigma_2... \sigma_n} = 2(m-1) + (2s+t(2s-2))(m-2) + t(m-s-1).$$

After simplifying the above, and adding $m$ to the result will produce the required expression for the size of the supersequence.
\end{proof}

For $s=1$, $s=2$, $s=3$ and $s=4$, these yield expressions $m^2 -2m + 4$, $m^2 - \frac {7}{3}m +\frac {19}{3}$, $m^2 - \frac {12}{5}m +\frac {38}{5}$ and $m^2 - \frac {17}{7}m +\frac {61}{7}$ respectively. The first and second of these are simply the Newey's and Radomirovic's expressions. Substitute $m=25$ into the last expression gives the result $573$, as demonstrated in Example 22 in the previous section.

The rest of the section is devoted to construct supersequence for $m$ outside of the congruence class $m \equiv 4 \pmod{2s-1}$. We need the following definition and Lemma from [8].

\begin{flushleft} \textbf{Definition 24.} Given any sequence $\sigma$ over a set $A$ and any letter $a$ from $A$ with $a \in \sigma$, then denote by $\sigma_{a...}$ the consecutive subsequence of $\sigma$ starting with the first occurence of $a$ in $\sigma$ until the end of $\sigma$, and with all occurrences of $a$ in $\sigma$ removed.
\end{flushleft} 

For example, if $\sigma = \langle 1, 2, 3, 4, 5, 4, 3, 2, 1 \rangle$, then the first occurence of $3$ is at the third position. Therefore $\sigma_{3...}$ denotes the subsequence starting at the third position until the end of $\sigma$, and with all occurrences of $3$ removed. So, we have $\sigma_{3...}= \langle 4, 5, 4, 2, 1 \rangle$.

\begin{flushleft} \textbf{Lemma 25.} Given a supersequence $\sigma$ over a set $A$ and any letter $a$ from $A$, then $\sigma_{a...}$ is a supersequence over $A \setminus \{a\}$.
\end{flushleft}
\begin{proof} Given any permutation $\rho$ of $A \setminus \{a\}$, then by supersequence nature of $\sigma$, we have $a \cdot \rho < \sigma$. Since  $\sigma_{a...}$ is a subsequence of $\sigma$ starting at the first occurence of $a$, we must have $\rho <\sigma_{a...}$. \end{proof}

We are ready to prove the following main results of the paper.

\begin{flushleft} \textbf{Theorem 26.} Given an integer $m \ge 4s+2$, there exists a supersequence of size at most $\lceil m^2 - \frac {5s-3}{2s-1}m +\frac {4s^3 -10s^2+21s-11}{2s-1} \rceil$ over a set of $m$ letters.
\end{flushleft}

\begin{proof} 

For any integer $m \ge 4s+2$, if $n=m-1$ satisfies the condition of Theorem 23, then we are done. Otherwise, pick the smallest integer $l$ greater than $m$ such that $l \equiv 4 \pmod{2s-1}$. So, we have $l-m \le 2s-2$. Let $\sigma$ be a supersequence over a set of $l$ letters as given by Theorem 23, constructed from $T_s(n)$ where $n=l-1$. We will repeatedly apply Lemma 25, at most $2s - 2$ many times, in order to construct a supersequence over a set of $m$ letters. Firstly, we pick $a$ to be the last non-skip letter to appear first time in $\sigma$ to create $\sigma_{a...}$. Then, we again pick $b$ to be the last non-skip letter to appear first time in $\sigma_{a...}$ to create $\sigma_{a...,b...}$. The process continues until we get a supersequence over a set of $m$ letters.

Observe that concatenating the sequences in $T_s(n)$ and removal of the skip letters results in a sequence that simply repeats $\langle 1,2,...,n-s+1 \rangle$. Recall $T_s(n)$ is defined from $n$-many sequences $\sigma_1, ... \sigma_n$. Since every non-skip letter is removed from at most one $\sigma_k$, so each non-skip letter occurs at least $n-1$ times in $\sigma$.

Thus, with $l=n+1$, when applying Lemma 25 the first time, $l-s$ non-skip letters including the initial $x$, appear before $a$. In addition, $a$ appears at least $l-2$ times. So the length of the supersequence will be reduced by at least $2l-(s+2)$. Both of these numbers reduce by one for each application of the Lemma, so at the $j$th step, the length reduces by at least $2l-(s+2j)$. Thus, if $m=l-d$, after summing $d$ terms and simplifying, we have a supersequence of size at most $\lceil m^2 - \frac {5s-3}{2s-1}m +\frac {2s^2+9s-7 +d(2s^2-4s+2)}{2s-1} \rceil$. This increases with $d$. Substituting the maximum value, $d=2s-2$ yields the result.
\end{proof}

The above proof only gives a near-optimal construction. An optimal approach will need to count the skip letters removed at each step $j$ as well. Due to the way that skip letters are shuffled, this will yield a piecewise function that is more complicated to optimize.

\begin{flushleft} \textbf{Corollary 27.} Given any real number $\epsilon > 0$, there exists a constant $C_{\epsilon}$ such that for all large enough $m$, there is a supersequence of length $\lceil m^2 - (\frac {5}{2}-\epsilon)m + C_{\epsilon} \rceil$ over a set of $m$ letters.
\end{flushleft}
\begin{proof} 
Since $\lim_{s \to \infty} \frac {5s-3}{2s-1} = \frac {5}{2}$, so for any $\epsilon > 0$, it is possible to find $s$ such that $ \frac {5s-3}{2s-1} > (\frac {5}{2}-\epsilon)$. Theorem 26 provides the supersequence, and $C_{\epsilon}$ is given by the third term of the expression.
\end{proof}

By choosing a reasonable $s$ for each $m$, it is possible to have the following bound.

\begin{flushleft} \textbf{Theorem 28.} For every $m$, there is a supersequence of size  $\lceil m^2 -  \frac {5}{2}m + \frac{3}{2}  (\frac{m}{2})^{\frac {2}{3}} +  (\frac{m}{2}) ^{\frac {1}{3}} + 7\rceil$.
\end{flushleft}
\begin{proof} 
A quick check shows that the result holds for $m < 10$. For each $m \ge 10$ , we have $m \ge 4s+2$, as required to apply Theorem 26. Let $s = \lceil \frac{1}{2} ((\frac{m}{2}) ^{\frac {1}{3}}+1) \rceil = \frac{1}{2} ((\frac{m}{2}) ^{\frac {1}{3}}+1) + \delta$, where $0 \le \delta < 1$. Rewrite the expression of Theorem 26 by polynomial division to obtain $ m^2 - \frac {5}{2}m + \frac {1}{2(2s-1)}m + 2s^2 -4s +\frac {17}{2} - \frac {5}{2(2s-1)} $. Substituting $s$ into the expression, we note that $ \frac {1}{2(2s-1)}m \le   (\frac{m}{2})^{\frac {2}{3}}$. After simplifying, we have a size that is at most 
$$ m^2 - \frac {5}{2}m + \frac{3}{2} (\frac{m}{2})^{\frac {2}{3}} - (\frac{m}{2})^{\frac {1}{3}} + 2\delta  (\frac{m}{2})^{\frac {1}{3}} - \frac{3}{2}-2\delta + 2\delta^2 +\frac {17}{2} $$
This gives the result as required.
\end{proof}

\begin{flushleft} \textbf{Acknowledgement:} The author thanks the referee for suggesting improvements on the results and the presentation.
\end{flushleft}

\begin{biog}
\item[Oliver Tan]
\begin{affil}

OliverTan333@gmail.com
\end{affil}

\end{biog}

\end{document}